\documentclass[11pt,a4paper]{article}

\usepackage{amsmath,amssymb}
\newcommand{\R}{\mathbb{R}}

\newcommand{\n}{\nabla}

\newcommand{\p}{\partial}

\newcommand{\g}{\int_{\R}}

\newcommand{\e}{\epsilon}
\newcommand{\va}{\varphi}

\newtheorem{theorem}{Theorem}

\newtheorem{remarka}{Remark}

\newtheorem{lemme}{Lemma}

\title{Weak-Strong uniqueness for compressible Navier-Stokes system with degenerate viscosity coefficient and  vacuum in one dimension}
\author{Boris Haspot 
\thanks{Ceremade UMR CNRS 7534
Universit\'e  Paris  Dauphine,
Place du MarŽchal DeLattre De Tassigny
75775 PARIS CEDEX 16 , haspot@ceremade.dauphine.fr }}
\date{}
\begin{document}
\maketitle




\begin{abstract}
We prove weak-strong uniqueness results for the compressible Navier-Stokes system with degenerate viscosity coefficient and with vacuum in one dimension. In other words, we give conditions on the weak solution constructed in \cite{Jiu} so that it is unique. The novelty consists in dealing with initial density $\rho_0$ which contains  vacuum. To do this we use the notion of relative entropy developed recently by Germain, Feireisl et al and Mellet and Vasseur (see \cite{PG,Fei,15}) combined with a new formulation of the compressible system (\cite{cras,CPAM,CPAM1,para}) (more precisely we introduce a new effective velocity which makes the system parabolic on the density and hyperbolic on this velocity).

\end{abstract}

\section{Introduction}
\subsection{Presentation of the equation}
We are interested in this paper by the Cauchy problem for the compressible Navier-Stokes equation with degenerate viscosity in the euclidian space $\R$. The system reads:
\begin{equation}
\begin{cases}
\begin{aligned}
&\p_t\rho+\p_x(\rho u)=0,\\
&\p_t(\rho u)+\p_x(\rho u^2)-\p_x(\mu (\rho)\p_x u)+\p_x P(\rho)=0.
\end{aligned}
\end{cases}
\label{1}
\end{equation}
with possibly degenerate viscosity coefficient. The unknowns $\rho$ and $u$ stand for the density and the velocity of the fluid, they are respectively $\R^+$ and $\R$ valued and they are functions of the space variable $x$ and of the time variable $t$. \\
Throughout the paper, we will assume that the pressure $P(\rho)$ obeys a $\gamma$ type law:
\begin{equation}
P(\rho)=\rho^\gamma,\;\gamma>1.
\end{equation}
Following the idea of \cite{cras,para,CPAM,CPAM1}, setting $v=u+\p_x \va(\rho)$ with $\va'(\rho)=\frac{\mu(\rho)}{\rho^2}$ we have (we refer to the appendix of \cite{global} for more details on the computations):
\begin{equation}
\begin{cases}
\begin{aligned}
&\p_t\rho-\p_{x}(\frac{\mu(\rho)}{\rho}\p_x\rho)+\p_x(\rho v)=0,\\
&\rho \p_t v+\rho u\p_x v+\p_x P(\rho)=0.
\end{aligned}
\end{cases}
\label{11}
\end{equation}
Let us mention that this change of unknown transform the system (\ref{1}) as a parabolic equation on the density and a transport equation for the velocity (we will see in the sequel that $v$ is not so far to verify a damped transport equation).It seems very surprising to observe that at the inverse of $u$ which has a parabolic behavior, $v$  has a hyperbolic behavior. In some sense in $1D$ the compressible Navier Stokes equations are a compressible Euler system with a viscous regularization on the density of the type $-\p_{x}(\frac{\mu(\rho)}{\rho}\p_x\rho)$\\
In the literature the authors are often interested in the constant viscosity case, however physically the viscosity of a gas depends on the temperature and  on the density (in the isentropic case). Let us mention the case of the Chapman-Enskog viscosity law (see \cite{CC70}) or the case of monoatomic gas ($\gamma=\frac{5}{3}$) where $\mu(\rho)=\rho^{\frac{1}{3}}$. More generally, $\mu(\rho)$ is expected to vanish as a power of the $\rho$ on the vacuum. In this paper we are going to deal with degenerate viscosity coefficients which can be written under the form $\mu(\rho)=\mu\rho^\alpha$ for some $\alpha>0$. In particular the case $\mu(\rho)=\mu$ and $\mu(\rho)=\mu\rho$ are included in our result (with $\mu$ positive constant).\\
\subsection{Weak and strong solutions in one dimension}

Let us start by recalling some results on the existence of solutions for the one dimension case when the viscosity is constant positive. The existence of global weak solutions was first obtained by Kazhikhov and Shelukin \cite{KS77} for smooth enough data close to the equilibrium (bounded away from zero). The case of discontinuous data (still bounded way from zero) was studied by Shelukin \cite{She82,She83,She84} and then by Serre \cite{Ser86a,Ser86b} and Hoff \cite{Hof87}. First results dealing with vacuum were also obtained by Shelukin  \cite{She86}. In \cite{Hof98}, Hoff extends the previous results by proving the existence of global weak solution with large discontinuous initial data having different limits at $x=\pm \infty$. In passing let us mention that the existence of global weak solution has been proved for the first time by P-L Lions in any dimension in \cite{14} and the result has been later refined by Feireisl et al \cite{FNP01}.\\
Concerning the existence of global weak solution with degenerate viscosity, the first results are due to Bresch and Desjardins in the multi d case where they introduce a new entropy modulo a algebraic relation between the viscosity coefficients (see \cite{BD}). It allows us to prove the stability of the global weak solution when we assume that the pressure is a cold pressure. This result has been extend in the framework of a $\gamma$ law by Mellet and Vasseur in \cite{MV06}. To do this the authors use an additional entropy giving a gain of integrability on the velocity which is sufficient in terms of compactness to deal with the terms $\rho u\otimes u$. Let us mention that all these results concern the multi d case $N\geq 2$ and that the problem of existence of global weak solution remains open. Indeed in \cite{MV06} the stability of the global weak solution is proved, however it seems difficult to construct global approximate solutions which verify uniformly all the different entropies. The existence of global strong solution for degenerate viscosity coefficient $\mu(\rho)=\rho^{\alpha}$ with $\alpha>\frac{1}{2}$) has been proved in one dimension by Jiu and Xin in \cite{Jiu}. They use the stability results proved in \cite{MV06} and they are able in the 1d case to construct global approximate solution of the system verifying uniformly all the entropies of \cite{MV06}. To do this they regularize the system by introducing a viscosity of the form $\mu_\delta(\rho)=\mu(\rho)+\delta\rho^{\theta}$ with $\theta<\frac{1}{2}$ (with such viscosity Mellet and Vasseur have proved the existence of global strong solution in \cite{15}).\\
\\
Concerning the uniqueness of the solution, Solonnikov  in \cite{Sol76} obtained the existence of strong solution for smooth initial data in finite time. However, the regularity may blow up when the solution approach from the vacuum. A natural question is to understand if the vacuum may arise in finite time. Hoff and Smoller (\cite{HS01}) show that any weak solution of the Navier-Stokes equation in one space dimension do not exhibit vacuum states, provided that no vacuum states are present initially. \\
Let us mention now some results  on the 1D compressible Navier-Stokes equations when the viscosity coefficient depends on the density. This case has been studied in
 \cite{LXY98}, \cite{OMNM02}, \cite{YYZ01}, \cite{YZ02}, \cite{JXZ} and \cite{LLX} in the framework of initial density admitting vacuum, more precisely the initial density is compactly supported and the authors are concerned by the evolution of the free boundary delimiting the vacuum. This is exactly the case that we are going to study in the sequel.\\
Let us mention the opposite case when the initial density is far away from the vacuum,  Mellet and Vasseur have proved in \cite{15} the existence of global strong solution with large initial data when $\mu(\rho)=\rho^{\alpha}$ with $0\leq\alpha<\frac{1}{2}$. This result has been extended in \cite{global} to the case of the shallow water equation ($\mu(\rho)=\mu\rho$).
\subsection{Weak-Strong uniqueness}
The idea of weak strong uniqueness is the following: assume that a weak solution has regular initial data such that it exists a strong solution in finite time associated to these initial data, then we can prove its uniqueness in the class of finite energy solutions. Weak strong uniqueness gives in particular conditions under which the equations is well behaved, and weak solutions are unique.\\
The first weak-strong uniqueness results were obtained by Prodi \cite{16} and Serrin \cite{17} for the incompressible Navier-Stokes equation. Many results of weak strong solution have been developed
these last years in different topic fluid mechanics, kinetic equation and a natural tool to deal with this type of problem is the so called relative entropy (see section \ref{section3} for a definition) . It seems that this notion has been first used to prove weak strong uniqueness results by Dafermos in \cite{4} who was considering conservation laws. Mellet and Vasseur in \cite{15} have proved using relative entropy weak strong uniqueness when the viscosity is non degenerate, it means $\mu(\rho)\geq\mu>0$. This relative entropy has also been used in the framework of kinetic equation, Berthelin and Vasseur in \cite{1} proved convergence of kinetic models to the isentropic compressible model.\\
The existence of weak-strong solution in the class of the weak solution has been proved for the first time by Germain in \cite{PG} using relative entropy in the multi d case. This result has been extended by Feireisl et al in \cite{Fei} by constructing suitable weak solution verifying the relative entropy inequality proposed in \cite{PG}.\\
In this paper we are interested in extending the result of \cite{15} to the case of degenerate viscosity coefficient with initial data admitting vacuum. To do this we are going to combined the notion of relative entropy and the new formulation of the compressible Navier-Stokes equation, see the system (\ref{11}). Up our knowledge it is the first result of weak strong uniqueness with degenerate viscosity coefficient and with vacuum for compressible Navier-Stokes equation, this is due to the fact that in one dimension the system (\ref{11}) is hyperbolic for the new unknown $v$. $v$ verifies a Euler equation.
\section{Main result}
We present below our main theorem of weak strong uniqueness for (\ref{11}). It will be a convenient short-hand to denote:
$$L^p_TL^q=L^p([0,T],L^q(\R)).$$
We denote in the sequel the energy ${\cal E}(\rho,u)$ and ${\cal E}(\rho,v)$ which is used in \cite{Jiu}:
$$
\begin{aligned}
&{\cal E}(\rho,u)(t)=\int_{\R}[\frac{1}{2}\rho(t,x)|u(t,x)|^2+\frac{1}{\gamma-1}\rho^{\gamma}(t,x)] dx+\int^t_0\int_{\R}\mu(\rho)|\p_x u|^2(s,x)ds dx,\\
&{\cal E}(\rho,v)(t)=\int_{\R}[\frac{1}{2}\rho(t,x)|v(t,x)|^2+\frac{1}{\gamma-1}\rho^{\gamma}(t,x)] dx+\int^t_0\int_{\R}\p_x P(\rho)\p_x (\frac{\mu(\rho)}{\rho^2}) ds dx.
\end{aligned}
$$
\begin{theorem}
Let $P(\rho)=a\rho^\gamma$ with $\gamma>1$ and $\mu(\rho)=\mu\rho^\gamma$. Assume that the initial data $(\rho_0,u_0)$ verify:
$$\rho_0\in L^1\cap L^\gamma,\;\sqrt{\rho_0}u_0\in L^2,\;\sqrt{\rho_0}v_0\in L^2\;\;\mbox{and}\;\;\rho_0^{\frac{1}{2+\delta}}u_0\in L^{2+\delta},$$
with $\delta>0$ arbitrarily small.
A solution $(\bar{\rho},\bar{v})$ of system (\ref{11}) is unique on $[0,T]$ in the set of solutions $(\rho,v)$ such that: 
\begin{equation}
\sup_{t\in(0,T)}\|\rho(t,\cdot)\|_{L^1}<+\infty,\;\sup_{t\in(0,T)}{\cal E}(\rho,u)(t)<+\infty\;\;\mbox{and}\;\;\sup_{t\in(0,T)}{\cal E}(\rho,v)(t)<+\infty,
\label{4}
\end{equation}
provided it satisfies $\bar{u}, \p_x\bar{u},\p_x\bar{v}\in L^1_T(L^\infty)$ and $\bar{u}\in L^{\infty}((0,T)\times\R)$.

\label{theo1}
\end{theorem}
\begin{remarka}
Let us mention that the assumption on $\bar{\rho}$ and $\bar{u}$ are verified in \cite{LLX,YZ02,Y1} in finite time in the case of a free boundary problem.
\end{remarka}
Our main result theorem \ref{theo1} is proved in section \ref{section3} where we recall the notion of relative entropy applied to the system (\ref{11}). It allows us to deal with the vacuum since the momentum equation in (\ref{11}) is only a Euler equation in one dimension. It allows us to avoid all the difficulty related to the vacuum when we have a degenerate viscosity coefficient. 
\section{Proof of theorem \ref{theo1}}
\label{section3}
\subsection{Relative entropy, Weak strong uniqueness}
Let us study the compressible Navier-Stokes system (\ref{11}) with
$v=u+\p_x \va(\rho)$ where $\va'(\rho)=\frac{\mu(\rho)}{\rho^2}$:
\begin{equation}
\begin{cases}
\begin{aligned}
&\p_t\rho-p_{x}(\frac{\mu(\rho)}{\rho}\p_x\rho)+\p_x(\rho v)=0,\\
&\rho \p_t v+\rho u\p_x v+\p_x P(\rho)=0.
\end{aligned}
\end{cases}
\label{11a}
\end{equation}
Let us consider two solution $(\rho,v)$ and $(\bar{\rho},\bar{v})$, we set:
$$
\begin{aligned}
&U=v-\bar{v},\;R=\rho-\bar{\rho}\;\;\mbox{and}\;\;U_1=u-\bar{u},\\
&F(\bar{\rho},R)=\frac{1}{\gamma}(R+\bar{\rho})^{\gamma}-\bar{\rho}^{\gamma-1}R-\frac{1}{\gamma}\bar{\rho}^\gamma.
\end{aligned}
$$
First subtracting the mass conservation of (\ref{1}) for $(\rho,u)$ and $(\bar{\rho},\bar{u})$ gives:
\begin{equation}
\p_t R+\p_x(\rho U_1)+\p_x(R\bar{u})=0.
\label{6}
\end{equation}
In a similar way, we can deal with the momentum equation and it yields:
$$
\begin{aligned}
&(\rho\p_t +\rho u\p_x)U+\p_x a\rho^{\gamma}-\p_x a\bar{\rho}^\gamma=-R(\p_t\bar{v}+\bar{u}\p_x\bar{v})-\rho U_1\p_x \bar{v}.
\end{aligned}
$$
Next using the fact that: 
$$\p_t\bar{v}+\bar{u}\p_x\bar{v}=-a\frac{\p_x\bar{\rho}^\gamma}{\bar{\rho}}.$$
we have:
\begin{equation}
\begin{aligned}
&(\rho\p_t +\rho u\p_x)U+\p_x a\rho^{\gamma}-a\frac{\rho}{\bar{\rho}}\p_x\bar{\rho}^\gamma=-\rho U_1\p_x \bar{v}.
\end{aligned}
\label{3}
\end{equation}
Next we have by multiplying the equation  (\ref{3}) by $U=U_1+\frac{\mu}{\gamma-1}\p_x(\rho^{\gamma-1}-\bar{\rho}^{\gamma-1})$. Following \cite{PG} we have:
$$
\begin{aligned}
&\int^t_0\g(\p_x a\rho^{\gamma}-a\frac{\rho}{\bar{\rho}}\p_x \bar{\rho}^\gamma).U_1 dx ds=\frac{\gamma}{\gamma-1}\g \rho U\cdot\p_x (\rho^{\gamma-1}-\bar{\rho}^{\gamma-1})dxds\\
&=-\frac{\gamma}{\gamma-1}\int^t_0 \g\p_x(\rho U)(\rho^{\gamma-1}-\bar{\rho}^{\gamma-1})dxds\;\;\mbox{by integrations by parts}\\
&=\frac{\gamma}{\gamma-1}\int^t_0 \g(\p_t R+\p_x(R\bar{u}))(\rho^{\gamma-1}-\bar{\rho}^{\gamma-1})dxds\;\;\mbox{by using (\ref{6})}\\
&=\frac{\gamma}{\gamma-1}\int^t_0 \g \p_t R\,\frac{\p}{\p R} F dxds+\int^t_0\g\bar{u}\p_x R\,\frac{\p}{\p R} F dxds+\int^t_0\g\p_x \bar{u} \,R\,\frac{\p}{\p R} F dxds,\\
&=\frac{\gamma}{\gamma-1}\int^t_0 \g \p_t F dxds-\g \p_t\bar{\rho} \,\frac{\p}{\p \bar{\rho}} F dxds+\int^t_0\g\bar{u}\p_x F dx-\int^t_0\g\bar{u}\p_x\bar{\rho}\frac{\p}{\p \bar{\rho}} F dx\\
&\hspace{9cm} +\int^t_0\g\p_x\bar{u} \,R\,\frac{\p}{\p R} F dxds,\\
&=\frac{\gamma}{\gamma-1}\int^t_0 \g \p_t F dxds+\int^t_0\g \p_x \bar{u}(-F+\bar{\rho}\frac{\p}{\p\bar{\rho}}F+R\,\frac{\p}{\p R} F) dxds,\\
&=\frac{\gamma}{\gamma-1}\int^t_0 \g \p_t F dxds+\gamma \int^t_0\g\p_x\bar{u}\,F dxds.
\end{aligned}
$$
\begin{remarka}
Let us observe that all the expressions written above converge properly, indeed let us deal with:
$$
\begin{aligned}
&\p_x a\rho^{\gamma}U_1=a\gamma \rho^{\gamma-\frac{3}{2}}\p_x\rho \,\sqrt{\rho}U_1.
\end{aligned}
$$
But we know that $\gamma \rho^{\gamma-\frac{3}{2}}\p_x\rho$ is in $L^\infty_T(L^2)$ via the energy ${\cal E}(\rho,v)$. And we have $\sqrt{\rho}u$ which belongs to $L^\infty_T(L^2)$. It is enough to conclude. Similarly since $\sqrt{\rho}$ is in $L^\infty_T(L^2)$ we show that $\gamma \rho^{\gamma-\frac{3}{2}}\p_x\rho \,\sqrt{\rho}\bar{u}$ is in $L^1_T(L^1)$ since $\bar{u}$ is in $L^1_T(L^\infty)$.\\
In the same way we have:
$$\frac{\rho}{\bar{\rho}}\p_x \bar{\rho}^\gamma \,U_1=\gamma\bar{\rho}^{\gamma-2}\p_x\bar{\rho}\sqrt{\rho}\sqrt{\rho}U.$$
Since $\sqrt{\rho}u$ is in $L^\infty_T(L^2)$, $\bar{\rho}^{\gamma-2}\p_x\bar{\rho}$ in $L^1_T(L\infty)$ and $\sqrt{\rho}$ in $L^\infty_T(L^2)$ we conclude by H\"older's inequality. We proceed similarly for the other term assuming that $\bar{u}$ belongs to $L^\infty_T(L^\infty)$.
\end{remarka}
Now we have to deal with the part $(\p_x a\rho^{\gamma}-a\frac{\rho}{\bar{\rho}}\p_x \bar{\rho}^\gamma)\frac{\mu}{\gamma-1}\p_x(\rho^{\gamma-1}-\bar{\rho}^{\gamma-1})$ with $\mu(\rho)=\mu\rho^m$ and $m=\gamma$:
$$
\begin{aligned}
&a\frac{\mu}{\gamma-1}\g(\p_x \rho^{\gamma}-\frac{\rho}{\bar{\rho}}\p_x \bar{\rho}^\gamma)\p_x(\rho^{\gamma-1}-\bar{\rho}^{\gamma-1})dx\\
&= a\frac{\mu}{\gamma-1}\gamma(\gamma-1)\g|\rho^{\gamma-\frac{3}{2}}\p_x\rho-\sqrt{\rho}\bar{\rho}^{\gamma-2}\p_x\bar{\rho}|^2 dx
\end{aligned}
$$
\begin{remarka}
Let us observe that in order to define the previous quantity we need to be sure that 
$\rho^{m-\frac{3}{2}}\p_x\rho$ and $\sqrt{\rho}\bar{\rho}^{m-2}\p_x\bar{\rho}$ are in $L^2_T(L^2(\R))$. Let us mention that this is the case via the energy inequality (\ref{4}) and the condition on $\bar{\rho}$.
\end{remarka}
Finally multiplying the momentum equation of (\ref{3}) by $U$ we obtain:
$$
\begin{aligned}
&\p_t\|\sqrt{\rho}U\|_{L^2}^2+\frac{\gamma}{\gamma-1}\p_t\|F(\bar{\rho},R)\|_{L^1}\\
&\hspace{3cm}+ a\mu\gamma
\|\rho^{m-\frac{3}{2}}\p_x\rho-\sqrt{\rho}\bar{\rho}^{m-2}\p_x\bar{\rho}\|_{L^2}^2\\
&\leq -\gamma \g\p_x\bar{u}\,F(\bar{\rho},R) dx-\g\rho U_1\cdot\p_x \bar{v}.U dx
\end{aligned}
$$
We have then:
$$
\begin{aligned}
&\g(\rho U_1\p_x \bar{v}).U dx=\g(\rho U\p_x \bar{v}).U dx-\mu \g \sqrt{\rho} \big((\rho^{\gamma-\frac{3}{2}}\p_x\rho-\sqrt{\rho}\bar{\rho}^{\gamma-2}\p_x\bar{\rho})\p_x \bar{v}\big).U dx
\end{aligned}
$$
By a basic Gronwall inequality and a bootstrap argument we conclude that if $\p_x\bar{u}$ and $\p_x\bar{v}$ belong to $L_T^1(L^\infty)$ then we have $\sqrt{\rho}v=\sqrt{\rho}\bar{v}$  and $\rho=\bar{\rho}$. $\blacksquare$


\end{document}